\documentclass[a4paper,draft,12pt,reqno]{article}
\usepackage{amsmath,amssymb,amsthm,xspace,a4wide,amscd,latexsym,amssymb}
\usepackage[all]{xy}
\usepackage{enumerate}

\theoremstyle{plain}
\newtheorem{theorem}{Theorem}[section]
\newtheorem{lemma}[theorem]{Lemma}

\newcommand{\HH}{{}^{\infty}_{\hspace{1mm}\lambda}\mathcal{H}_{\mu}^1}

\begin{document} 

\title{Representation type of ${}^{\infty}_{\hspace{2mm}\lambda}\mathcal{H}_{\mu}^1$}
\author{Yuriy Drozd and Volodymyr Mazorchuk}
\date{}
\maketitle

\begin{abstract}
For a semi-simple finite-dimensional complex Lie algebra $\mathfrak{g}$ we classify the 
representation type of the associative algebras associated with the categories $\HH$ of 
Harish-Chandra bimodules for $\mathfrak{g}$.
\end{abstract}

\section{The result}\label{s1}

Let $\mathfrak{g}$ be a simple finite-dimensional complex Lie algebra
with a fixed triangular decoposition, $\mathfrak{g}=\mathfrak{n}_-\oplus
\mathfrak{h}\oplus\mathfrak{n}_+$, let $\lambda$ and $\mu$ be two
dominant and integral (but not necessarily regular) weights, let 
$U(\mathfrak{g})$ be the universal enveloping algebra of $\mathfrak{g}$,
and let $Z(\mathfrak{g})$ be the center of $U(\mathfrak{g})$. 
Denote by  $\chi_{\lambda}$ and $\chi_{\mu}$ the central characters of
the Verma modules $\Delta(\lambda)$ and $\Delta(\mu)$ respectively.
Let further $\HH$ denote the full subcategory of the category of
all $U(\mathfrak{g})$-bimodules, which consists of all $X$ satisfying
the following conditions (see \cite[Kapitel~6]{Ja}):
\begin{enumerate}[(1)]
\item $X$ is finitely generated as a bimodule;
\item $X$ is algebraic, that is $X$ is a direct sum of finite-dimensional
$\mathfrak{g}$-modules with respect to the diagonal action
$g\mapsto(g,\sigma(g))$, where $\sigma$ is the Chevalley involution on
$\mathfrak{g}$;
\item $x(z-\chi_{\mu}(z))=0$ for all $x\in X$ and $z\in Z(\mathfrak{g})$;
\item for every $x\in X$ and $z\in Z(\mathfrak{g})$ there exists
$k\in\mathbb{N}$ such that $(z-\chi_{\lambda}(z))^kx=0$.
\end{enumerate}
For regular $\mu$ the category $\HH$ is equivalent to a block of the 
BGG category $\mathcal{O}$, associated with the triangular decomposition above, 
see \cite{BGG}. For singular $\mu$ the category $\HH$ is equivalent to 
a block of the parabolic generalization $\mathcal{O}(\mathfrak{p},\Lambda)$
of $\mathcal{O}$, studied in  \cite{FKM}. Moreover, from \cite{FKM,So1} it follows
that every block of $\mathcal{O}$ and $\mathcal{O}(\mathfrak{p},\Lambda)$ is
equivalent to some $\HH$. Every $\HH$ is equivalent to the module category of a
properly stratified finite-dimensional associative algebra. The regular blocks 
of $\HH$ can be used to categorify a parabolic Hecke module,  see \cite{MS}.

Let $\mathbf{W}$ be the Weyl group of $\mathfrak{g}$ and  $\rho$ be the half 
of the sum of all positive roots of $\mathfrak{g}$. Then $\mathbf{W}$ acts on
$\mathfrak{h}^*$ in the usual way and we recall the following {\em dot-action} 
of $\mathbf{W}$ on $\mathfrak{h}^*$: $w\cdot \nu=w(\nu+\rho)-\rho$.
Let $\mathbf{G}\subset \mathbf{W}$ be the stabilizer of $\lambda$
with respect to the dot-action, and $\mathbf{H}\subset \mathbf{W}$ be the 
stabilizer of $\mu$ with respect to the dot-action. We will say that the triple
$(\mathbf{W},\mathbf{G},\mathbf{H})$ is associated to $\HH$. In the present 
paper we classify the categories $\HH$ according to their representation type 
in terms of the associated triples, thus extending the results of \cite{FPN,BKM,GP}. 
Let $(\mathbf{W},\mathbf{G},\mathbf{H})$ be the triple, associated to $\HH$, and
$(\mathbf{W},\mathbf{G}',\mathbf{H}')$ be the triple, associated to some
${}^{\infty}_{\hspace{0.7mm}\lambda'}\mathcal{H}_{\mu'}^1$. Then
from  \cite[Theorem~5.9]{BG} and \cite[Theorem~11]{So1} it follows that 
$\HH$ and ${}^{\infty}_{\hspace{0.7mm}\lambda'}\mathcal{H}_{\mu'}^1$  are equivalent 
if there exists an automorphism, $\varphi$, of the Coxeter system $(\mathbf{W},S)$,
where $S$ is the set of simple reflections associated to our triangular decomposition, 
such that $\varphi(\mathbf{G})=\mathbf{G}'$ and $\varphi(\mathbf{H})=\mathbf{H}'$.
By the {\em Coxeter type} of a triple, $(\mathbf{W},\mathbf{G},\mathbf{H})$, we mean the
triple that consists of the Coxeter types of the corresponding components of
$(\mathbf{W},\mathbf{G},\mathbf{H})$. Note that, in general, the Coxeter type of the 
triple does not determine the triple in a unique way (for example, one can  compare 
the cases \eqref{tm.1.5}, \eqref{tm.2.4} and \eqref{tm.2.5} in the formulation
of Theorem~\ref{tmain} below). Our main result is the following statement:

\begin{theorem}\label{tmain}
\begin{enumerate}[(1)]
\item\label{tm.1} The category $\HH$ is of finite type
if and only if the Coxeter type of the associated triple is
\begin{enumerate}[(a)]
\item\label{tm.1.1} any and $\mathbf{W}=\mathbf{G}$;
\item\label{tm.1.2} $(A_n,A_{n-1},A_n)$, $(B_n,B_{n-1},B_n)$, $(C_n,C_{n-1},C_n)$,  or 
$(G_2,A_1,G_2)$;
\item\label{tm.1.3} $(A_1,e,e)$;
\item\label{tm.1.4} $(A_n,A_{n-1},A_{n-1})$;
\item\label{tm.1.5} $(A_n,A_{n-1},A_{n-2})$, where $A_{n-2}$ is obtained from $A_{n}$ 
by taking away the first and the last roots;
\item\label{tm.1.6} $(B_2,A_{1},A_{1})$ or $(C_2,A_1,A_1)$, and 
$\mathbf{G}=\mathbf{H}$ (in both cases);
\item\label{tm.1.7} $(B_n,B_{n-1},B_{n-1})$ or $(C_n,C_{n-1},C_{n-1})$, where $n\geq 3$;
\item\label{tm.1.8} $(A_2,A_{1},e)$.
\end{enumerate}
\item\label{tm.2} The category $\HH$ is tame if and only if the Coxeter type of 
the associated triple is
\begin{enumerate}[(a)]
\item\label{tm.2.1} $(A_3,A_1\times A_1,A_3)$, $(A_2,e,A_2)$, $(B_2,e,B_2)$, 
$(G_2,e,G_2)$, $(B_3,A_2,B_3)$, $(C_3,A_2,C_3)$, or $(D_n,D_{n-1},D_n)$ where $n\geq 4$;
\item\label{tm.2.2} $(B_2,A_1,A_1)$ or $(C_2,A_1,A_1)$, and $\mathbf{G}\neq\mathbf{H}$
(in both cases);
\item\label{tm.2.3} $(A_n,A_{n-1},A_{1}\times A_{n-2})$, $n>2$;
\item\label{tm.2.4} $(A_n,A_{n-1},A_{n-2})$, $n>2$, where $A_{n-2}$ is included into 
$A_{n-1}$ and contains either the first or the last root of $A_n$;
\item\label{tm.2.5} $(A_n,A_{n-1},A_{n-2})$, $n>2$, where $A_{n-2}$ is not included into 
$A_{n-1}$;
\item\label{tm.2.7} $(A_3,A_2,e)$, $(B_2,A_1,e)$, $(C_2,A_1,e)$.
\end{enumerate}
\item\label{tm.3} In all other cases the category $\HH$ is wild.
\end{enumerate}
\end{theorem}

For regular $\mu$ Theorem~\ref{tmain} gives the classification of the 
representation type of the blocks of the category $\mathcal{O}$ obtained in 
\cite{FPN} (see also \cite{BKM} for a different proof). Formally, we do not use any 
results from \cite{FPN} and \cite{BKM}, however, the main idea of our proof is
similar to the one of \cite{BKM}.

In the case $\mathbf{H}=\mathbf{W}$ (i.e. $\mu$ is most singular) Theorem~\ref{tmain} 
reduces to the classification of the representation type for the algebra
$\mathtt{C}(\mathbf{W},\mathbf{G})$ of $\mathbf{G}$-invariants in the coinvariant 
algebra associated to $\mathbf{W}$. This result was obtained in \cite{GP} and, in fact, 
our argument in the present paper is based upon it. 

The last important ingredient in the proof of Theorem~\ref{tmain}, the latter being 
presented in  Section~\ref{s3}, is the classification of the representation type of all 
centralizer subalgebras in the Auslander algebra $\mathtt{A}_n$ of $\Bbbk[x]/(x^n)$. 
This classification is given in Section~\ref{s2}. Two series of centralizer subalgebras, 
namely those considered in Lemma~\ref{lh7} and Lemma~\ref{lh8}, seem to be rather 
interesting and non-trivial. 

The paper finishes with an extension of Theorem~\ref{tmain} to the case of a semi-simple 
Lie algebra $\mathfrak{g}$. This is presented in Section~\ref{s4}, where one more 
interesting tame algebra arises.

We would like to finish the introduction with a remark that just recently a first step towards 
the classification of the representation type of the blocks of Rocha-Caridi's parabolic 
analogue $\mathcal{O}_S$ of $\mathcal{O}$ was made in \cite{BN}. The next step would be to 
complete this classification and then to classify the representation type of the ``mixed'' 
version of $\mathcal{O}_S$ and $\mathcal{O}(\mathfrak{p},\Lambda)$. As the results of 
\cite{BN} and of the present paper suggest, this might give some interesting tame algebras 
in a natural way.

\section{Representation type of the centralizer subalgebras in the Auslander
algebra of $\Bbbk[x]/(x^n)$}\label{s2}

In the paper we will compose arrows of the quiver algebras from the right to the 
left. Let $\Bbbk$ be an algebraically closed field. Recall that, according to 
\cite{Dr3}, every finite-dimensional associative $\Bbbk$-algebra has either finite, 
tame or wild representation type. In what follows we will call the latter statement 
the {\em Tame and Wild Theorem}. The algebras, which are not of finite representation 
type, are said to be of infinite representation type.

Let $A=(A_{ob},A_{mor})$ be a $\Bbbk$-linear category. An $A$-module, $M$, is a functor
from $A$ to the category of $\Bbbk$-vector spaces. In particular, for $x\in A_{ob}$ and
$\alpha\in A_{mor}$ we will denote by $M(x)$ and $M(\alpha)$ the images of $x$ and $\alpha$
under $M$ respectively.

For a positive integer $n>1$ let $\mathtt{A}_n$ be the algebra given by the 
following quiver with relations:
\begin{displaymath}
\xymatrix{
1\ar@/^/[rr]^{a_1} && 2\ar@/^/[rr]^{a_2}\ar@/^/[ll]^{b_1} && \dots
\ar@/^/[ll]^{b_2}\ar@/^/[rr]^{a_{n-1}} && n\ar@/^/[ll]^{b_{n-1}} 
}\quad\quad
\begin{array}{ll}
a_{i}b_{i}=b_{i+1}a_{i+1}, & i=1,\dots,n-2,\\
a_{n-1}b_{n-1}=0.
\end{array}
\end{displaymath}
The algebra $\mathtt{A}_n$ is the Auslander algebra of $\Bbbk[x]/(x^n)$
(see for example \cite[Section~7]{DR}).
For $X\subset \{2,3,\dots,n\}$ let $e_X$ denote the direct sum of all primitive
idempotents of $\mathtt{A}_n$, which corrrespond to the vertexes from $\{1\}\cup X$.
Set $\mathtt{A}_n^{X}=e_X\mathtt{A}_n e_X$. The main result of this section is the
following:

\begin{theorem}\label{taus}
\begin{enumerate}[(i)]
\item\label{taus.1} The algebra $\mathtt{A}_n^{X}$ has finite representation type
if and only if $X\subset \{2,n\}$.
\item\label{taus.2} The algebra $\mathtt{A}_n^{X}$ has tame representation type if 
and only if either $n>3$ and $X=\{3\}$, $\{2,3\}$, $\{n-1\}$, $\{n-1,n\}$, or 
$n=4$ and $X=\{2,3,4\}$.
\item\label{taus.3} The algebra $\mathtt{A}_n^{X}$ is wild in all other cases.
\end{enumerate}
\end{theorem}

To prove Theorem~\ref{taus} we will need the following lemmas:

\begin{lemma}\label{lh1}
The algebra $\mathtt{A}^{\{m\}}_n$ has infinite representation type for
$m\in\{3,\dots,n-1\}$ and $n\geq 4$.
\end{lemma}

\begin{proof}
The algebra $\mathtt{A}^{\{m\}}_n$ is given by the following quiver with relations:
\begin{equation}\label{eqh1}
\xymatrix{ 
1\ar@/^/[rr]^a\ar@(ul,dl)[]_{x}   && m\ar@/^/[ll]^b\ar@(ur,dr)[]^{y}
} \quad\quad\quad 
\begin{array}{ll}
ax=ya, & xb=by,\\
ab=y^{m-1}, & ba=x^{m-1},\\
y^{n-m+1}=0, 
\end{array}
\end{equation}
where $x=b_1a_1$, $y=b_ma_m$, $a=a_{m-1}\dots a_1$, $b=b_1\dots b_{m-1}$. Modulo the
square of the radical $\mathtt{A}^{\{m\}}_n$ gives rise to the following diagram of infinite type:
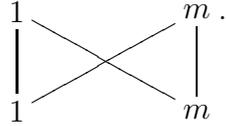
\begin{displaymath}
\xymatrix{
 1\ar@{-}[rrd]\ar@{-}[d] && m\ar@{-}[d] \\
 1\ar@{-}[rru] && m
}.
\end{displaymath}
Hence $\mathtt{A}^{\{m\}}_n$ has infinite representation type as well.
\end{proof}

\begin{lemma}\label{lh2}
The algebra $\mathtt{A}^X_n$ is wild for $X=\{3,m\}$, where $m>4$.
\end{lemma}

\begin{proof}
In this case the algebra $\mathtt{A}^X_n$ is given by the following quiver with relations:
\begin{displaymath}
\xymatrix{ 
1\ar@/^/[rr]^a\ar@(ul,dl)[]_{x}   && 3 \ar@/^/[ll]^b\ar@/^/[rr]^s\ar@(ul,ur)[]^{y} 
&& m\ar@/^/[ll]^t\ar@(ur,dr)[]^{z}
} \quad\quad 
\begin{array}{ll}
ax=ya, & xb=by,\\
sy=zs, & yt=tz,\\
ab=y^2,& ba=x^2,\\
st=z^{m-3},&ts=y^{m-3},\\
z^{n-m+1}=0, 
\end{array}
\end{displaymath}
where $x=b_1a_1$, $y=b_3a_3$, $z=b_ma_m$, $a=a_2a_1$, $b=b_1b_2$,
$s=a_{m-1}\dots a_3$, $t=b_3\dots b_{m-1}$. Note that $z=0$ if $m=n$.
Modulo the square of the radical $\mathtt{A}^{X}_n$ gives rise to the following diagram:
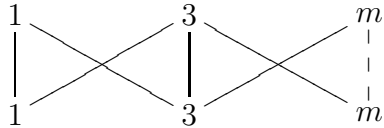
\begin{displaymath}
\xymatrix{
 1\ar@{-}[rrd]\ar@{-}[d] && 3\ar@{-}[d]\ar@{-}[rrd] && m\ar@{--}[d] \\
 1\ar@{-}[rru] && 3\ar@{-}[rru] && m
}
\end{displaymath}
(where the dashed line disappears in the case $m=n$). With or without the dashed
line the diagram is not an extended Dynkin quiver and hence is wild
(see \cite{DF,DR0}). Hence $\mathtt{A}^{X}_n$ is wild as well.
\end{proof}

\begin{lemma}\label{lh3}
The algebra $\mathtt{A}^X_n$ is wild for $X=\{2,n-1\}$ and $n\geq 5$.
\end{lemma}

\begin{proof}
To make the quivers in the proof below look better we set $m=n-1$.
The algebra $\mathtt{A}^X_n$ is given by the following quiver with relations:
\begin{displaymath}
\xymatrix{ 
1\ar@/^/[rr]^a   && 2 \ar@/^/[ll]^b\ar@/^/[rr]^s 
&& m\ar@/^/[ll]^t\ar@(ur,dr)[]^{x}
} \quad\quad 
\begin{array}{ll}
sab=xs, & abt=tx,\\
st=0,  &ts=(ab)^{n-3},\\
x^{2}=0, 
\end{array}
\end{displaymath}
where $a=a_1$, $b=b_1$, $s=a_{n-2}\dots a_2$, $t=b_2\dots b_{n-2}$, $x=b_{n-1}a_{n-1}$. 
The universal covering of $\mathtt{A}^{X}_n$ has the wild fragment (a hereditary algebra,
whose underlined quiver is not an extended Dynkin diagram, see \cite{DF,DR0}) indicated by the
dotted arrows in the following picture:
\begin{displaymath}
\xymatrix{
\dots  && \dots &&&& \dots &&&& \\
 1\ar[rr]^{a} && 2 \ar[rrrr]^{s}\ar[lld]_{b}&&&& m \ar@{.>}[d]^{x}
 \ar@{.>}[dddllll]|->>>>>>>>>>{t}\\
 1\ar@{.>}[rr]^{a} && 2 \ar@{.>}[rrrr]^{s}\ar@{.>}[lld]_{b}&&&& m 
 \ar[d]^{x}\ar[dddllll]|->>>>>>>>>>{t}\\
 1\ar@{.>}[rr]^{a} && 2 \ar[rrrr]^{s}\ar@{.>}[lld]_{b}&&&& m \ar[d]^{x}\\
 1\ar@{.>}[rr]^{a} && 2 \ar[rrrr]^{s}\ar[lld]_{b}&&&& m \ar[d]^{x}\\
 1\ar[rr]^{a} && 2 \ar[rrrr]^{s}&&&& m \\
\dots && \dots &&&& \dots\\
}
\end{displaymath}
Hence $\mathtt{A}^{X}_n$ is wild as well.
\end{proof}

\begin{lemma}\label{lh4}
The algebra $\mathtt{A}^{\{3,4\}}_5$ is wild.
\end{lemma}

\begin{proof}
The algebra $\mathtt{A}^{\{3,4\}}_5$ is given by the following quiver with relations:
\begin{displaymath}
\xymatrix{ 
1\ar@/^/[rr]^a \ar@(ul,dl)[]^{x}  && 3 \ar@/^/[ll]^b\ar@/^/[rr]^s 
&& 4\ar@/^/[ll]^t
} \quad\quad 
\begin{array}{ll}
ax=tsa, & xb=bts,\\
ba=x^2,&ab=(ts)^{2},\\
(st)^{2}=0,
\end{array}
\end{displaymath}
where $a=a_2a_1$, $b=b_1b_2$, $s=a_3$, $t=b_3$, $x=b_{1}a_{1}$. 
The universal covering of $\mathtt{A}^{\{3,4\}}_5$ has the wild fragment (a hereditary algebra,
whose underlined quiver is not an extended Dynkin diagram, see \cite{DF,DR0}) indicated by the
dotted arrows in the following picture:
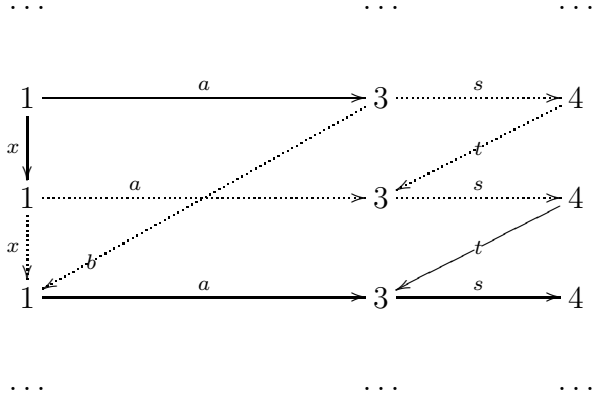
\begin{displaymath}
\xymatrix{
\dots  &&&& \dots && \dots &&&& \\
1\ar[rrrr]^{a}\ar[d]_x &&&& 3\ar@{.>}[rr]^{s}\ar@{.>}[ddllll]|->>>>>>>{b} 
&& 4\ar@{.>}[dll]|-t\\
1\ar@{.>}[rrrr]^>>>>>>>>>>>>>>>>>>>>>>>>>>>>>>{a}\ar@{.>}[d]_x 
&&&& 3\ar@{.>}[rr]^{s} && 4\ar[dll]|-t\\
1\ar[rrrr]^{a} &&&& 3\ar[rr]^{s} && 4\\
\dots &&&& \dots && \dots\\
}
\end{displaymath}
Hence $\mathtt{A}^{\{3,4\}}_5$ is wild as well.
\end{proof}

\begin{lemma}\label{lh6}
The algebra $\mathtt{A}^{\{m\}}_n$ is wild for $m\in\{4,\dots,n-2\}$
and $n\geq 6$.
\end{lemma}

\begin{proof}
The algebra $\mathtt{A}^{\{m\}}_n$ is given by \eqref{eqh1}. We consider its quotient
$\mathtt{B}$ given by the additional relations $x^3=y^3=ab=ba=0$ (which is possible because of
our restrictions on $m$ and $n$). Then the universal covering of $\mathtt{B}$ exists
and has the following fragment,
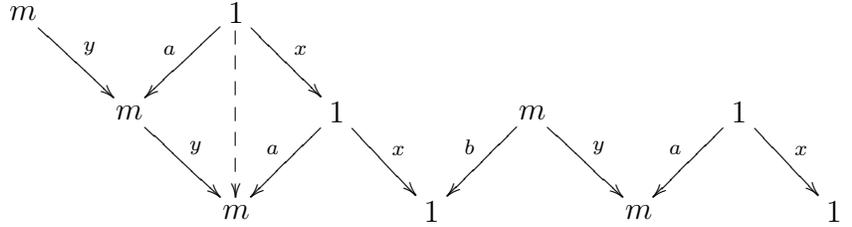
\begin{displaymath}
\xymatrix{
m\ar[dr]^y && 1\ar[dr]^x\ar[dl]_a\ar@{-->}[dd] \\
&\ar[dr]^y m && 1\ar[dr]^x\ar[dl]_a && m\ar[dr]^y\ar[dl]_b && 1\ar[dr]^x\ar[dl]_a \\
&& m && 1 && m && 1 
},
\end{displaymath}
which is wild by \cite{Un}.
This implies that $\mathtt{B}$ and hence $\mathtt{A}^{\{m\}}_n$ is wild.
\end{proof}

\begin{lemma}\label{lh7}
The algebra $\mathtt{A}^{\{2,n\}}_n$, $n\geq 2$, 
is  of finite representation  type.
\end{lemma}

\begin{proof}
For $n=2,3$ the statement follows from \cite[Section~7]{DR}.
The algebra $\mathtt{A}^{\{2,n\}}_n$, $n\geq 4$, is given by the 
following quiver with relations:
\begin{equation}\label{equequ}  
\xymatrix{ 1 \ar@/^/[r]^a & 2 \ar@/^/[l]^b \ar@/^/[r]^u & n \ar@/^/[l]^v  }
\qquad  uv=uab=abv=0,\ vu=(ab)^{n-2}, 
\end{equation} 
where $a=a_1$, $b=b_1$, $u=a_{n-1}\dots a_2$, $v=b_2\dots b_{n-1}$. Note that these 
relations imply $(ab)^{n-1}=(ba)^n=0$. The projective $\mathtt{A}^{\{2,n\}}_n$-module 
$P(1)$ is injective, so we can replace $\mathtt{A}^{\{2,n\}}_n$ by 
$\mathtt{A}'=\mathtt{A}^{\{2,n\}}_n/\mathrm{soc}(P(1))= \mathtt{A}^{\{2,n\}}_n/((ba)^{n-1})$, 
which has the same indecomposable modules except $P(1)$, see \cite[Lemma~9.2.2]{DK}. So from 
now on we consider the algebra $\mathtt{A}'$, i.e. add the relation $(ba)^{n-1}=0$ to 
\eqref{equequ}. The algebra $\mathtt{A}'$ has a simply connected covering 
$\widetilde{\mathtt{A}}$, see \cite{BoGa}, which is the category, given by the following 
quiver with relations (we show the case $n=5$, in the general case the arrow starting at 
$n_k$ ends at $2_{n-2+k}$):
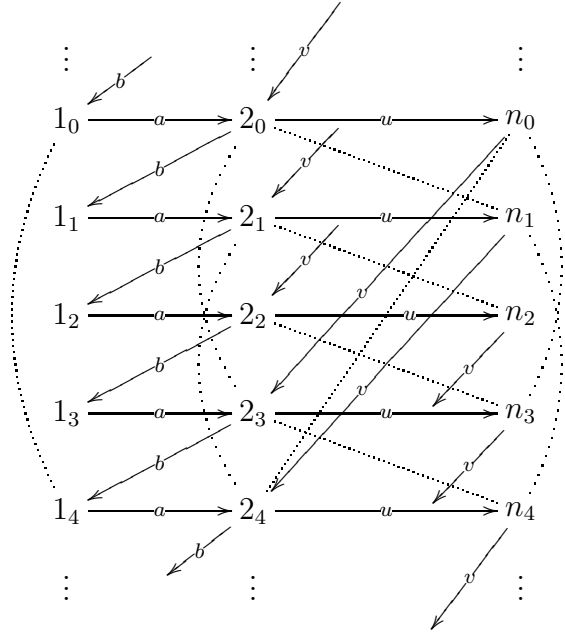
\begin{displaymath}
\xymatrix@R=1.5ex{  
 &&&\ar[ddl]|-v&&      \\
\vdots & \ar[dl]|-b &\vdots&&&\vdots \\
1_0 \ar@/_4ex/@{.}[dddddddd] \ar[rr]|-a && 2_0 \ar[rrr]|-u \ar[ddll]|-b \ar@/_4ex/@{.}[dddddd]
\ar@{.}[ddrrr] &\ar[ddl]|-v&& n_0 \ar[ddddddlll]|->>>>>>>>>>>>>>>>>v 
\ar@/^3ex/@{.}[dddddd] \ar@{.}[ddddddddlll]\\  
&&&& &&\\
1_1 \ar[rr]|-a && 2_1 \ar[rrr]|-u \ar[ddll]|-b \ar@{.}[ddrrr] \ar@/_4ex/@{.}[dddddd]
&\ar[ddl]|-v&& n_1 \ar[ddddddlll]|->>>>>>>>>>>>>>>>>v \ar@/^3ex/@{.}[dddddd] \\  
&&&&&& \\
1_2 \ar[rr]|-a && 2_2 \ar[rrr]|->>>>>>>>>>>u \ar[ddll]|-b \ar@{.}[ddrrr]&&& n_2 \ar[ddl]|-v \\  
&&&  &&&\\
1_3 \ar[rr]|-a && 2_3 \ar[rrr]|-u \ar[ddll]|-b \ar@{.}[ddrrr]&&& n_3 \ar[ddl]|-v \\  
&&&  &&&\\
1_4 \ar[rr]|-a && 2_4 \ar[rrr]|-u \ar[dl]|-b &&& n_4 \ar[ddl]|-v \\  
\vdots&&\vdots&&&\vdots      \\
&&&&&      \\
}
\end{displaymath} 
We omit the indices at the arrows $a,b,u,v$. They satisfy the same relations 
as in $\mathtt{A}'$, which are shown by the dotted lines. Consider the full subcategory 
$\mathtt{B}_m$ of $\widetilde{\mathtt{A}}$ with the set of objects 
$\mathbf{S}=\{1_k, m\le k\leq m+n-1;\,\,2_k, m\leq k\leq m+n-2;\,\,n_m\}$.
Let $M$ be an $\widetilde{\mathtt{A}}$-module, $N_m$ be its restriction to 
$\mathtt{B}_m$, $N_m=\bigoplus_{i=1}^sK_i$, where $K_i$ are indecomposable 
$\mathtt{B}_m$-modules. It is well known that every $K_i$ is completely determined by 
the subset of objects $\mathbf{S}_i=\{x\,|\,K_i(x)\ne0\}$ and if $1_m\in\mathbf{S}_i$, 
then $1_{m+n-1}\notin\mathbf{S}_i$. Moreover, all $K_i(x)$ with $x\in\mathbf{S}_i$ are 
one-dimensional and all arrows between these objects correspond to the identity maps.  
Since $uab=abv=0$, $K_i$ splits out of the whole module $M$ whenever 
$\mathbf{S}_i\supseteq\{2_m,2_{m+n-2}\}$. Suppose that, for every integer $m$, $N_m$ does 
not contain such direct summands. It implies that $M(vu)=0$. Therefore $M$ can be 
considered as a module over $\overline{\mathtt{A}}$, where $\overline{\mathtt{A}}$ is given 
by the following quiver
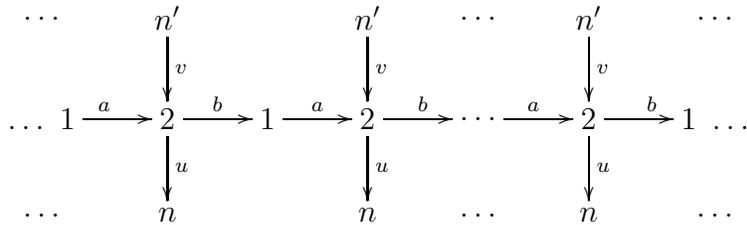
\begin{displaymath} 
\xymatrix{ \dots & n' \ar[d]^v && n' \ar[d]^v & \dots & n' \ar[d]^v &\dots\\
\dots\ 1 \ar[r]^a & 2 \ar[r]^b\ar[d]^u & 1 \ar[r]^a & 2 \ar[r]^b \ar[d]^u
&\dots \ar[r]^a & 2\ar[r]^b \ar[d]^u & 1\ \dots \\
\dots & n && n & \dots & n &\dots }
\end{displaymath} 
with relations $ uv=uab=abv=(ab)^{n-2}=0$. One easily checks that any indecomposable
representation of $\overline{\mathtt{A}}$ is at most of dimension $2n-5$. Hence, 
$\overline{\mathtt{A}}$ is representation (locally) finite, i.e. for every object 
$x\in\overline{\mathtt{A}}$ there are only finitely many indecomposable
representations $M$ with $M(x)\ne0$. By \cite{BoGa}, the algebra  $\mathtt{A}^{\{2,n\}}_n$
is representation (locally) finite as well, which completes the proof.
\end{proof}

\begin{lemma}\label{lh8}
The algebra $\mathtt{A}^{\{n-1,n\}}_n$, $n>3$, is tame.
\end{lemma}

\begin{proof}
For $q=n-1$ the algebra $\mathtt{A}_n^{\{q,n\}}$ is given by the 
following quiver with relations
\begin{displaymath}
\xymatrix{ 1 \ar@(ul,dl)[]_c \ar@/^/[r]^u & {q} \ar@/^/[l]^v \ar@/^/[r]^a &
n \ar@/^/[l]^b }\quad\quad  c^n=ab=uv=0,\ vu=c^{n-2},\ cv=vba,\ uc=bau,
\end{displaymath} 
where $c=b_1a_1,\,a=a_{q},\,b=b_{q},\,u=a_{n-2}\dots a_1,\,v=b_1\dots b_{n-2}$.
The projective module $P(1)$ is also injective, hence, using \cite[Lemma~9.2.2]{DK} as
it was done in the proof of Lemma~\ref{lh7}, we can replace  $\mathtt{A}$ 
by $\mathtt{A}'=\mathtt{A}/\mathrm{soc}(P(1))=\mathtt{A}/(c^{q})$. Let $M$ be an 
$\mathtt{A}'$-module. Choose a basis in $M(1)$ so that the matrix  $C=M(c)$ is in the 
Jordan normal form, or, further,
\begin{displaymath} 
M(c)=\bigoplus_{i=1}^{q} J_i\otimes I_{m_i},  
\end{displaymath} 
where $J_i$ is the nilpotent Jordan block of size $i\times i$ and $I_{m_i}$ is the identity
matrix of size $m_i\times m_i$ (here $m_i$ is just the number of Jordan blocks of size $i$). 
Thus
\begin{displaymath} 
J_i\otimes I_m =  \begin{pmatrix}
0 & I_m & 0 & \dots & 0 & 0 \\ 0& 0& I_m & \dots & 0 & 0 \\ \hdotsfor 6 \\
0&0&0&\dots &0&I_m \\ 0&0&0&\dots &0&0
\end{pmatrix}_{i\times i}
\end{displaymath}
(here $i\times i$ means $i$ boxes times $i$ boxes, each of size $m_i$). Choose bases 
in $M(q)$ and $M(n)$ such that the matrices $A=M(a)$ and $B=M(b)$ are of the form
\begin{displaymath} 
A= \begin{pmatrix}
0 & 0 & 0 & I & 0 \\ 0 & 0 & 0 & 0 & I \\ 0&0&0&0&0 \\ 0&0&0&0&0
\end{pmatrix}, \qquad
B= \begin{pmatrix}
0 & I & 0 & 0 \\0 & 0 & 0 & I \\ 0&0&0&0 \\ 0&0&0&0 \\ 0&0&0&0
\end{pmatrix},
\end{displaymath} 
where the vertical (horizontal) stripes of $A$ are of the same size as the horizontal
(respectively, vertical) stripes of $B$, and $I$ is the identity matrix; 
we do not specify these sizes here. Set $r=nq/2$; it is the number of the horizontal 
and vertical stripes in $C$. Then $M(u)$
and $M(v)$ can be considered as block matrices: $M(u)=U=(U_k^{ij})_{5\times r}$
and $M(v)=V=(V_{ij}^k)_{r\times5}$,
where $k=1,\dots,5$ correspond to the $k$-th horizontal stripe of $B$; 
$i=1,\dots,q,\ j=1,\dots,i$, and the stripe $(ij)$ corresponds to the $j$-th horizontal
stripe of the matrix $J_i\otimes I_{m_i}$ in the decomposition of $C$. The conditions
$uc=bau$ and $cv=vba$ imply that for $i>1$ the only nonzero blocks $U_k^{ij}$ and 
$V_{ij}^k$ can be
\begin{align*} 
U_k^{ii}\quad &\text{ and }\ U_1^{i,i-1}=U_5^{ii} \\
V^k_{i1}\quad &\text{ and }\ V^5_{i2}=V^1_{i1}.
\end{align*} 
Moreover, we also have $U_5^{11}=V^1_{11}=0$. Changing bases in the spaces $M(x)$,
$x=1,q,n$, so that the matrices $A,B,C$ remain of the same form, we can
replace $U$ and $V$ respectively by $T^{-1}US$ and $S^{-1}VT$, where $S,T$ are
invertible matrices of the appropriate sizes such that $SA=AS$ and $TU=UQ,\,QV=VT$
for an invertible matrix $Q$. We also consider $S$ and $T$ as block matrices:
$S=(S_{st}^{ij})_{r\times r}$ and $T=(T_l^k)_{5\times5}$ with respect to the division
of $A,B,C$. Then the conditions above can be rewritten as follows: 
\begin{itemize}
\item 	$S_{st}^{ij}$ can only be nonzero if $i-j<s-t$ or $i-j=s-t,\,s\le i$;
\item  $S_{st}^{ij}=S_{st'}^{ij'}$ if $t-j=t'-j'$;
\item  $T$ is block triangular: $T_l^k=0$ if $k<l$, and $T^1_1=T^5_5$;
\item  all diagonal blocks $S^{ij}_{ij}$ and $T_k^k$ are invertible.
\end{itemize} 
Especially, for the vertical stripes $U^{ii}$ and for the horizontal stripes $U_k$
of the matrix $U$ the following transformations are allowed:
\begin{enumerate}
\item\label{ob1}  Replace $U^{ii}$ by $U^{ii}Z$.
\item\label{ob2}  Replace $U_k$ by $ZU_k$, where $k=2,3,4$.
\item\label{ob3}  Replace $U_1$ and $U_5$ respectively by $ZU_1$ and $ZU_5$.
\item\label{ob4}  Replace $U^{ii}$ by $U^{ii}+U^{jj}Z$, where $j<i$.
\item\label{ob5}  Replace $U_k$ by $U_k+U_lZ$, where $k<l$.
\end{enumerate}
Here $Z$  denotes an arbitrary matrix of the appropriate size, moreover, in the 
cases \ref{ob1}--\ref{ob3}
it must be invertible. One can easily see that, using these transformations,
one can subdivide all blocks $U^{ii}_k$ into subblocks so that each stripe
contains at most one nonzero block, which is an identity matrix. Note that the
sizes of the horizontal substripes of $U_1$ and $U_5$ must be the same. 
Let $\Lambda^{ii}$ and $\Lambda_k$ be respectively the sets of the vertical and the
horizontal stripes of these subdivisions. Note that all stripes $U^{ij}$ must
be subdivided respectively to the subdivision of $U^{ii}$ and recall that
$U_1^{i,i-1}=U_5^{ii}$. Especially, there is a one-to-one correspondence
$\lambda\mapsto \lambda'$ between $\Lambda_5$ and $\Lambda_1$.

We make the respective subdivision of the blocks of the matrix $V$, too. 
The condition $UV=0$ implies that, whenever the $\lambda$-th vertical stripe of $U$
is nonzero ($\lambda\in\Lambda^{ii}$), the $\lambda$-th horizontal stripe of $V$ is zero. The
conditions $VU=C^{q}$ can be rewritten as 
\begin{displaymath} 
  V_{ij}U^{st}= \begin{cases}
  I &\text{ if }\ (i,j,s,t)=(q,1,q,q),\\
 0 &\text{ otherwise}.
 \end{cases}  
\end{displaymath} 
It implies that there are no zero vertical stripes in the new subdivision of $U^{q,q}$.
Moreover, if $\lambda\in\Lambda^{ii}$, $\mu\in\Lambda_k$, and the block $V^\lambda_\mu$ is 
nonzero, then the $\mu$-th vertical stripe of $U$ is zero if $i\ne q$; if $i=q$ this 
stripe contains exactly one non-zero block, namely, $U^\mu_\lambda=I$. We denote by 
$\overline{\Lambda}^{ii}$ and $\overline{\Lambda}_k$ the set of those stripes from 
$\Lambda^{ii}$ and $\Lambda_k$, which are not completely defines by these rules. Let
$\lambda\in\Lambda_5$, $\lambda'$ be the corresponding element of $\Lambda_1$. If the blocks
$U_\lambda^\mu$ and $U_{\lambda'}^{\mu'}$ are both nonzero, write $\mu\sim\mu'$. Note that 
there is at most one element $\mu'$ such that it holds, and $\mu'\ne\mu$.

One can verify that the sets $\overline{\Lambda}^{ii}$ and $\overline{\Lambda}_k$ can 
be linearly ordered so that, applying the transformations of the types \ref{ob1}--\ref{ob5} 
from above, we can 
replace a stripe $V^\lambda$ by $V^\lambda+V^{\lambda'} Z$ with $\lambda'<\lambda$ and a 
stripe $V_\mu$ by $V_\mu+ZV_{\mu'}$, where $\lambda'<\lambda,\,\mu'<\mu$ for any matrix $Z$ 
(of the appropriate size). We can also replace $V^\lambda$ by $V^\lambda Z$, where $Z$ is 
invertible, and replace simultaneously $V_\mu$ and $V_{\mu'}$, where $\mu'\sim\mu$, by 
$ZV_\mu$ and $ZV_{\mu'}$ (if $\mu'$ does not exist, just replace $V_\mu$ by $ZV_\mu$) with 
invertible $Z$. Therefore, we obtain a special sort of the matrix problems considered in 
\cite{Bo}, which is known to be tame. Hence, the algebra $\mathtt{A}_n^{\{q,n\}}$ is tame 
as well.
\end{proof}

\begin{proof}[Proof of Theorem~\ref{taus}.]
Lemma~\ref{lh7} and Lemma~\ref{lh1} imply Theorem~\ref{taus}\eqref{taus.1}. 
The statement of Theorem~\ref{taus}\eqref{taus.3} follows from
Theorem~\ref{taus}\eqref{taus.1} and Theorem~\ref{taus}\eqref{taus.2} using the
Tame and Wild Theorem.
Hence we have to prove Theorem~\ref{taus}\eqref{taus.2} only. 

It is known, see for example \cite{DR}, that $\mathtt{A}_n$ has finite representation 
type for $n\leq 3$, is tame for $n=4$, and is wild for all other $n$. This, in particular, 
proves Theorem~\ref{taus}\eqref{taus.2} for $n\leq 4$.

If $n\geq 6$ then from Lemma~\ref{lh6} it follows that if $\mathtt{A}^X_n$ is tame then
$X\subset\{2,3,n-1,n\}$. From Theorem~\ref{taus}\eqref{taus.1} we know that
$X\not\subset \{2,n\}$. From Lemma~\ref{lh2} it follows that $\{3,n-1\}\not\subset X$
and $\{3,n\}\not\subset X$. From Lemma~\ref{lh3} it follows that $\{2,n-1\}\not\subset X$.
This leaves us the cases $X=\{n-1,n\}$, $\{n-1\}$, $\{2,3\}$ and $\{3\}$.
In the first two cases $\mathtt{A}^X_n$ is tame by Lemma~\ref{lh8}. The algebra
$\mathtt{A}^{\{2,3\}}_n$, $n\geq 3$,  is given by the following quiver with relations:
\begin{displaymath}
\xymatrix{
1\ar@/^/[rr]^a && 2\ar@/^/[rr]^s\ar@/^/[ll]^b && 3\ar@/^/[ll]^t
}\quad\quad
\begin{array}{l}
ab=ts,\\
(st)^{n-2}=0,
\end{array}
\end{displaymath}
where $a=a_1$, $b=b_1$, $s=a_2$, $t=b_2$. 
For $n\geq 5$ this algebra is tame as a quotient of the classical tame problem from 
\cite{NR}. Hence $\mathtt{A}^{\{3\}}_n$ is tame as well. 

For $n=5$ Lemma~\ref{lh4} implies that $\mathtt{A}^X_n$ is wild if
$X\supset\{3,4\}$, Lemma~\ref{lh2} implies that $\mathtt{A}^X_n$ is 
wild if $X\supset\{3,5\}$, and Lemma~\ref{lh3} implies that $\mathtt{A}^X_n$ is 
wild if $X\supset\{2,4\}$. Above we have already shown that the algebras 
$\mathtt{A}^{\{2,3\}}_5$ is tame, and hence $\mathtt{A}^{\{3\}}_5$ is tame as well. 
Finally, that the algebras $\mathtt{A}^{\{4,5\}}_5$ and $\mathtt{A}^{\{4\}}_5$ are tame
follows from Lemma~\ref{lh8}. This completes the proof.
\end{proof}

\section{Proof of Theorem~\ref{tmain}}\label{s3}

We briefly recall the structure of $\HH$. We refer the reader to \cite{BG,So1,FKM,KM}
for details. By \cite[Theorem~5.9]{BG},
the category ${}^{\infty}_{\hspace{1mm}\lambda}\mathcal{H}_{0}^1$ is equivalent
to the block $\mathcal{O}_{\lambda}$ of the BGG category $\mathcal{O}$, \cite{BGG}. 
Let $\mathtt{O}(\mathbf{W},\mathbf{G})$ denote the basic associative algebra,
whose module category is equivalent to $\mathcal{O}_{\lambda}$.
The simple modules in $\mathcal{O}_{\lambda}$ are in 
natural bijection with the cosets $\mathbf{W}/\mathbf{G}$ (under this bijection the
coset $\mathbf{G}$ corresponds to the dominant highest weight). For 
$w\in \mathbf{W}$ let $L(w)$ denote the corresponding simple module in 
$\mathcal{O}_{\lambda}$, $P(w)$ be the projective cover of $L(w)$, $\Delta(w)$ be
the corresponding Verma module, and $I(w)$ be the injective envelope of $L(w)$. 
Then \cite{So1} implies that for the longest element $w_0\in\mathbf{W}$ one has
$\mathrm{End}_{\mathcal{O}_{\lambda}}(P(w_0))\cong\mathtt{C}(\mathbf{W},\mathbf{G})$
(recall that this is the subalgebra of $\mathbf{G}$-invariants in the coinvariant 
algebra, associated to $\mathbf{W}$). The left multiplication in $\mathbf{W}$ induces 
an action of $\mathbf{H}$ on the set $\mathbf{W}\cdot \lambda$. Let $P(\lambda,\mathbf{H})$ 
denote the direct sum of indecomposable projective modules that correspond to the longest 
elements in all orbits of this action. The category $\HH$ is equivalent, by \cite{KM}, 
to the module category over $\mathtt{B}(\mathbf{G},\mathbf{H})=
\mathrm{End}_{\mathcal{O}_{\lambda}}(P(\lambda,\mathbf{H}))$. From \cite{So1} it follows
that $\mathtt{B}(\mathbf{G},\mathbf{H})$ depends on $\mathbf{G}$ rather than on $\lambda$. 

We start with Theorem~\ref{tmain}\eqref{tm.1}, that is with the case of finite
representation type.

Note that $P(w_0)$ is always a direct summand of $P(\lambda,\mathbf{H})$. Hence
$\mathtt{C}(\mathbf{W},\mathbf{G})$ is a centralizer subalgebra of 
$\mathtt{B}(\mathbf{G},\mathbf{H})$. In particular, for $\HH$ to be of finite 
representation type, $\mathtt{C}(\mathbf{W},\mathbf{G})$ must be of finite representation 
type as well. According to \cite[Theorem~7.2]{GP}, $\mathtt{C}(\mathbf{W},\mathbf{G})$ is 
of finite representation type in the following cases:
\begin{enumerate}[(I)]
\item\label{s3.1} $\mathbf{W}=\mathbf{G}$;
\item\label{s3.2} $\mathbf{W}$ is of type $A_n$ and $\mathbf{G}$ is of type $A_{n-1}$;
\item\label{s3.3} $\mathbf{W}$ is of type $B_n$ and $\mathbf{G}$ is of type $B_{n-1}$;
\item\label{s3.4} $\mathbf{W}$ is of type $C_n$ and $\mathbf{G}$ is of type $C_{n-1}$;
\item\label{s3.5} $\mathbf{W}$ is of type $G_2$ and $\mathbf{G}$ is of type $A_1$.
\end{enumerate}
Moreover, in all these cases $\mathtt{C}(\mathbf{W},\mathbf{G})\cong\mathbb{C}[x]/(x^r)$,
where $r=[\mathbf{W}:\mathbf{G}]$. The last observation and \cite[Theorem~1]{FKM}
imply that in all the above cases the category $\mathcal{O}_{\lambda}$ is equivalent to 
$\mathtt{A}_r\mathrm{-mod}$. In particular, the algebra 
$\mathtt{B}(\mathbf{G},\mathbf{H})$ is isomorphic to $\mathtt{A}_r^X$ for appropriate $X$,
and, in the notation of Section~\ref{s2}, the algebra $\mathtt{C}(\mathbf{W},\mathbf{G})$
is the centralizer subalgebra, which corresponds to the vertex $1$.

The case \eqref{s3.1} gives Theorem~\ref{tmain}\eqref{tm.1.1}. 
In the cases \eqref{s3.2}, \eqref{s3.3}, \eqref{s3.4}, and \eqref{s3.5} it follows from 
Theorem~\ref{taus}\eqref{taus.1} that we have the following possibilities
for $\mathtt{B}(\mathbf{G},\mathbf{H})$:

{\bf $\mathtt{B}(\mathbf{G},\mathbf{H})$ has one simple module.} This implies
$\mathbf{W}=\mathbf{H}$ and gives Theorem~\ref{tmain}\eqref{tm.1.2}.

{\bf $\mathtt{B}(\mathbf{G},\mathbf{H})$ has two simple modules.}
These simples correspond either to the dominant and the anti-dominant weights 
in $\mathcal{O}_{\lambda}$ or to the anti-dominant weight and its neighbor. 
By a direct calculation we get the following: the case
$r=2$ gives Theorem~\ref{tmain}\eqref{tm.1.3}, and the case
$r>2$ gives Theorem~\ref{tmain}\eqref{tm.1.4}. 

{\bf $\mathtt{B}(\mathbf{G},\mathbf{H})$ has three simple modules.}
These simples correspond to the following weights in $\mathcal{O}_{\lambda}$:
the anti-dominant one, its neighbor, and the dominant one.
By a direct calculation we get the following: the case
$r=3$ gives Theorem~\ref{tmain}\eqref{tm.1.8}, and the case
$r>3$ gives Theorem~\ref{tmain}\eqref{tm.1.5}, Theorem~\ref{tmain}\eqref{tm.1.6},
and Theorem~\ref{tmain}\eqref{tm.1.7}. This proves Theorem~\ref{tmain}\eqref{tm.1}.

Let us now proceed with the tame case, that is with Theorem~\ref{tmain}\eqref{tm.2}. 
If $\mathtt{C}(\mathbf{W},\mathbf{G})$ is of finite representation type, that is in the
cases \eqref{s3.1}--\eqref{s3.5}, Theorem~\ref{taus}\eqref{taus.2} give us
the following possibilities for $\mathtt{B}(\mathbf{G},\mathbf{H})$:

{\bf $\mathtt{B}(\mathbf{G},\mathbf{H})$ has two simple modules.} 
These simples correspond to the following weights in $\mathcal{O}_{\lambda}$:
either the anti-dominant one and the neighbor of its neighbor, or
the anti-dominant one and the neighbor of the dominant one.
By a direct calculation we get that these cases lead to
Theorem~\ref{tmain}\eqref{tm.2.2} and Theorem~\ref{tmain}\eqref{tm.2.3}.

{\bf $\mathtt{B}(\mathbf{G},\mathbf{H})$ has three simple modules.} 
These simples correspond to the following weights in $\mathcal{O}_{\lambda}$:
either the anti-dominant one, its neighbor, and the neighbor of its neighbor, 
or the anti-dominant, its  neighbor and the dominant one.
By a direct calculation we get that these cases lead to
Theorem~\ref{tmain}\eqref{tm.2.4} and Theorem~\ref{tmain}\eqref{tm.2.5}.

{\bf $\mathtt{B}(\mathbf{G},\mathbf{H})$ has four simple modules.} In this
case $r=4$ and a direct calculation gives Theorem~\ref{tmain}\eqref{tm.2.7}.

The rest (that is Theorem~\ref{tmain}\eqref{tm.2.1}) should correspond to the case when
$\mathtt{C}(\mathbf{W},\mathbf{G})$ is tame. According to \cite[Theorem~7.2]{GP},
$\mathtt{C}(\mathbf{W},\mathbf{G})$ is tame in the following cases:
\begin{enumerate}[(I)]
\setcounter{enumi}{5}
\item\label{s3.n1} $\mathbf{W}$ has rank $2$ and $\mathbf{G}=\{e\}$;
\item\label{s3.n2} $\mathbf{W}$ is of type $A_3$ and $\mathbf{G}$ is of type $A_{1}\times A_1$;
\item\label{s3.n3} $\mathbf{W}$ is of type $B_3$ and $\mathbf{G}$ is of type $A_2$;
\item\label{s3.n4} $\mathbf{W}$ is of type $C_3$ and $\mathbf{G}$ is of type $A_{2}$;
\item\label{s3.n5} $\mathbf{W}$ is of type $D_n$ and $\mathbf{G}$ is of type $D_{n-1}$.
\end{enumerate}
For $\mathbf{W}=\mathbf{H}$ the cases \eqref{s3.n1}, \eqref{s3.n2}, \eqref{s3.n3}, 
\eqref{s3.n4}, and \eqref{s3.n5} give exactly Theorem~\ref{tmain}\eqref{tm.2.1}. Let
us now show that the rest is wild.

If $\mathbf{W}\neq\mathbf{H}$ then $\HH$ has at least two non-isomorphic indecomposable 
projective modules, one of which is $P(w_0)$ and the other one is some $P(w)$. 
We first consider the cases \eqref{s3.n2}, \eqref{s3.n3}, 
\eqref{s3.n4}, and \eqref{s3.n5}. In all these cases the restriction of the Bruhat order
to $\mathbf{W}/\mathbf{G}$ gives the following poset:
\begin{equation}\label{eqeqptm}
\xymatrix{
 & & & & u_1\ar@{-}[rd]& & &&\\
w_0 \ar@{-}[r]& w_1\ar@{-}[r] & \dots\ar@{-}[r] & w_s\ar@{-}[ru]\ar@{-}[rd] & & 
v_s\ar@{-}[r] & \dots\ar@{-}[r]& v_1\ar@{-}[r] & v_0 \\
 & & & & u_2\ar@{-}[ru]& & &&\\
}
\end{equation}
From \cite[Theorem~7.3]{GP} it follows that in all these cases 
the algebra $\mathtt{C}(\mathbf{W},\mathbf{G})$ has two  generators.

We consider the centralizer subalgebra 
$\mathtt{D}(w)=\mathrm{End}_{\mathcal{O}_{\lambda}}(P(w_0)\oplus P(w))$ and 
let $Q(w)$ denote the quotient of $\mathtt{D}(w)$ modulo the square of the
radical. Recall that the algebra $\mathtt{O}(\mathbf{W},\mathbf{G})$ is Koszul,
see \cite{BGS}, and hence the category $\mathcal{O}_{\lambda}$ is positively 
(Koszul) graded, see also \cite{St}. Hence $\mathtt{D}(w)$ is positively graded as well. 
We are going to show that $\mathtt{D}(w)$ is always wild. We start with the 
following statement.

\begin{lemma}\label{lnewmul}
Let $w\in\{w_1,\dots,w_s,u_1,u_2,v_0,\dots,v_s\}$. Then 
\begin{displaymath}
[P(v_0):L(w)]=
\begin{cases}
1, & w\in \{u_1,u_2,v_0,\dots,v_s,w_0\};\\
2, & w\in \{w_1,w_3,\dots,w_s\},
\end{cases}
\end{displaymath}
where $[P(v_0):L(w)]$ denotes the composition multiplicity.
\end{lemma}

\begin{proof}
By \cite{BGS} the category $\mathcal{O}_{\lambda}$ is Koszul dual to the regular 
block of the corresponding parabolic category of Rocha-Caridi, see \cite{RC}. Hence 
the multiplicity question for $\mathcal{O}_{\lambda}$ reduces, via the Koszul duality,
to the computation of the extensions in the parabolic case. The latter are given by 
Kazhdan-Lusztig polynomials and for the algebras of type \eqref{s3.n2}, \eqref{s3.n3}, 
\eqref{s3.n4}, and \eqref{s3.n5} these multiplicities are computed in \cite[\S~14]{ES}. 
The statement of our lemma follows directly from \cite[\S~14]{ES}.
\end{proof}

Since $L(w_0)$ is a simple Verma module, it occurs exactly one time in the composition 
series of $\Delta(w)$, which gives rise to a morphism, $\alpha:P(w_0)\to P(w)$. This 
morphism has the minimal possible degree (with respect to our positive grading) and hence
does not belong to the square of the radical. Further, the unique (now by the BGG 
reciprocity) occurrence of $\Delta(w)$ in the Verma flag of $P(w_0)$ gives a morphism, 
$\beta:P(w)\to P(w_0)$, which does not belong to the square of the radical either since it
again has the minimal possible degree. Now we will have to consider several cases.

{\bf Case~A.} Assume first that $w\in\{v_0,v_1,\dots,v_s\}$. The quiver of $Q(w)$ contains 
the arrows, corresponding to $\alpha$ and $\beta$. Moreover $Q(w)$ also contains two loops 
at the point $w_0$ which correspond to the generators of $\mathtt{C}(\mathbf{W},\mathbf{G})$. 
Passing, if necessary, to a quotient of $Q(w)$, we obtain the following configuration:
\begin{equation}\label{wconf2}
\xymatrix{
w_0\ar@/^/@{-}[d]\ar@/_/@{-}[d]\ar@{-}[rrd] && w\ar@{-}[lld]\\
w_0 && w
}.
\end{equation}
Since the underlined diagram is not an extended Dynkin diagram, the configuration
is wild, see \cite{DF,DR0}. This implies that $\mathtt{D}(w)$ and hence $\HH$ is wild in this case.

{\bf Case~B.} Consider now the case $w=u_1$ (the case $w=u_2$ is analogous).
Lemma~\ref{lnewmul} implies that in this case the multiplicity of
$L(w)$ in $\Delta(v_0)$ is $1$. Hence from \cite[Proposition~2.12]{Ba} it follows 
that $P(w)$ has simple socle $L(w_0)$, in particular, $P(w)$ is a submodule of 
$P(w_0)=I(w_0)$. Injectivity of $P(w_0)$ thus gives a surjection from 
$\mathrm{End}_{\mathcal{O}_{\lambda}}(P(w_0))\cong
\mathtt{C}(\mathbf{W},\mathbf{G})$ to $\mathrm{End}_{\mathcal{O}_{\lambda}}(P(w))$.
Note that, by \cite{So1}, $\mathrm{End}_{\mathcal{O}_{\lambda}}(P(w_0))$ is the center 
of $\mathtt{O}(\mathbf{W},\mathbf{G})$ and hence is central in 
$\mathtt{B}(\mathbf{G},\mathbf{H})$. We still have the elements $\alpha$ and $\beta$ as above,
which do not belong to the square of the radical. Further, using the embedding 
$P(w)\hookrightarrow P(w_0)$ one also obtains that $\alpha$ generates 
$\mathrm{Hom}_{\mathcal{O}_{\lambda}}(P(w_0),P(w))$ as a
$\mathtt{C}(\mathbf{W},\mathbf{G})$-module and $\beta$ generates 
$\mathrm{Hom}_{\mathcal{O}_{\lambda}}(P(w),P(w_0))$ as a
$\mathtt{C}(\mathbf{W},\mathbf{G})$-module.

With this notation, $\mathtt{D}(w)$ has the following quiver:
\begin{displaymath}
\xymatrix{
w_0\ar@/^/[rr]^{\alpha}\ar@(lu,ld)[]_{x} && w\ar@/^/[ll]^{\beta}\ar@(ru,rd)[]^{y}
}.
\end{displaymath}
Note that $\alpha$ is surjective as a homomorphism from
$\mathrm{End}_{\mathtt{D}(w)}(P(w_0))$ to $\mathrm{End}_{\mathtt{D}(w)}(P(w))$
since $P(w)$ has simple socle.
This and the fact that $\mathrm{End}_{\mathcal{O}_{\lambda}}(P(w_0))$ is central 
implies the relations $\alpha x=y\alpha$ and $\beta y=x\beta$. Using
\cite[7.12-7.16]{GP} one also easily gets the following additional relations: 
$y^{s+2}=0$, $\alpha\beta=c y^{s+1}$ for some $0\neq c\in\mathbb{C}$, $x\beta\alpha=
\beta\alpha x=0$ and $(\beta\alpha)^2=x^{2s+3}$. This implies that the universal 
covering of $\mathtt{D}(w)$ has the following fragment (shown for $s=1$):
\begin{equation}\label{eqeqnlm}
\xymatrix{
{\bf w_0}\ar[d]_{x}\ar[rr]^{\alpha}\ar@{--}[drr] && 
w\ar[d]_{y}\ar[ddll]|->>>>>>{\beta} \\
w_0\ar[rr]|->>>>{\alpha}\ar[d]_{x} && w \\
w_0 
}
\end{equation}
(here the dashed arrow indicates the commutativity of the corresponding square).
Evaluating the Tits form of this fragment at the point $(1,2,2,2,2)$, where 
$1$ is placed in the bold vertex, we obtain $-1<0$ implying that the fragment
\eqref{eqeqnlm} is wild (see for example \cite{CB,Dr}). Hence $\mathtt{D}(w)$ 
is wild as well.

{\bf Case~C.} Assume now that $w=w_i$, $i=2,\dots,s-1$. Hence, by  Lemma~\ref{lnewmul} 
the multiplicity of $L(w)$ in $P(v_0)$ is $2$. We will need the following lemma:

\begin{lemma}\label{hlem}
Let $A$ be a basic associative algebra, let $e$ be an idempotent of $A$ and
$f$ be a primitive direct summand of $e$. Assume that there exist two non-isomorphic
$A$-modules $M$ and $N$ satisfying the following properties:
\begin{enumerate}[(1)]
\item\label{hlem.1} both $M$ and $N$ have simple top and simple socles isomorphic to the
simple $A$-module $L^{A}(f)$, corresponding to $f$;
\item \label{hlem.2} $e\,\mathrm{rad}(M)/\mathrm{soc}(M)=e\,\mathrm{rad}(N)/\mathrm{soc}(N)=0$.
\end{enumerate}
Then $\dim \mathrm{Ext}_{eAe}^1(L^{eAe}(f),L^{eAe}(f))>1$.
\end{lemma}

\begin{proof}
Recall from \cite[Chapter~5]{Au} that $eAe\mathrm{-mod}$ is equivalent to the full
subcategory $\mathcal{M}$ of $A\mathrm{-mod}$, consisting of all $Ae$ approximations of 
modules from $A\mathrm{-mod}$. Let $M'$ and $N'$ be the $Ae$-approximations of $M$ and $N$ 
respectively. Both $M'$ and $N'$ are indecomposable since $M$ and $N$ are indecomposable
by \eqref{hlem.1}. Then the $eAe$-modules $eM'$ and $eN'$ are indecomposable as well, and,
because of \eqref{hlem.1} and \eqref{hlem.2}, both $eM'$ and $eN'$ have length two with 
both composition subquotients isomorphic to the simple $eAe$-module $L^{eAe}(f)$. 

Assume that $eM'\cong eN'$. Then, by \cite[Chapter~5]{Au}, any $eAe$-isomorphism between 
$eM'$ and $eN'$ induces an $A$-isomorphism between $M'$ and $N'$. From \eqref{hlem.1} we 
also have that the canonical maps $N\to N'$ and $M\to M'$ are injective, that is  we have
\begin{displaymath}
N\hookrightarrow N'\cong M'\hookleftarrow M.
\end{displaymath}
From \eqref{hlem.1}, the definition of the $Ae$-approximation, and the fact that 
$f$ is a direct summand of $e$, it follows that the image of $N$ in $N'$ coincides with the
trace of the projective module $Af$ in $N'$. Analogously the image of $M$ in $M'$ coincides 
with the trace of the projective module $Af$ in $M'$. This implies $M\cong N$, a
contradiction. The statement follows.
\end{proof}

Since we are not in the multiplicity-free case, from the Kazhdan-Lusztig Theorem 
it follows that the quiver of $\mathtt{O}(\mathbf{W},\mathbf{G})$ contains  more arrows 
than what is indicated on the diagram \eqref{eqeqptm}. Namely, from the results of 
\cite[\S~14]{ES} we have $\mathrm{Ext}_{\mathcal{O}_{\lambda}}^1(L(w),L(v_{i-1}))\neq 0$.
Note that $\mathrm{Ext}_{\mathcal{O}_{\lambda}}^1(L(w),L(w_{i+1}))\neq 0$
also follows from the Kazhdan-Lusztig Theorem since $w_i$ and $w_{i+1}$ are neighbors
(it follows from \cite[\S~14]{ES} as well).
Let now $u\in\{v_{i-1},w_{i+1}\}$. Then we can fix a non-zero element from
$\mathrm{Ext}_{\mathcal{O}_{\lambda}}^1(L(w),L(u))$. This means that $L(u)$ occurs in 
degree $1$ in the projective module $P(w)$. The module $P(w)$ has a Verma flag, and the above 
occurrence of $L(u)$ gives rise to an occurrence of $\Delta(u)$ as a subquotient of $P(w)$. 
Since $L(u)$ is in degree $1$ and $\mathcal{O}_{\lambda}$ is positively graded, we can
factor all the Verma subquotients of $P(w)$ except $\Delta(w)$ and $\Delta(u)$ out 
obtaining a non-split extension, $N(u)$ say, of $\Delta(u)$ by $\Delta(w)$. By duality, 
we have $\mathrm{Ext}_{\mathcal{O}_{\lambda}}^1(L(u),L(w))\neq 0$ as well, and, as $w<u$,
the module $L(w)$ occurs in degree $2$ in the module $N(u)$. This occurrence gives rise to a 
map from $N$ to the injective module $I(w)$. Let $N'(u)$ denote the image of this map. By 
construction, the module $N'(u)$ is an indecomposable module of Loewy length $3$ with
simple top and simple socle isomorphic to $L(w)$. Moreover,
$\mathrm{Rad}(N'(u))/\mathrm{Soc}(N'(u))$ (the latter is considered as an object
of $\mathcal{O}_{\lambda}$) does not contain $L(w)$ as a subquotient because of the 
quasi-hereditary vanishing $\mathrm{Ext}_{\mathcal{O}_{\lambda}}^1(L(w),L(w))=0$. Since 
$w\neq w_0,w_1$, all occurrences of $L(w_0)$ in $P(w)$ are in degrees $\geq 2$.  Hence
$\mathrm{Rad}(N'(u))/\mathrm{Soc}(N'(u))$ does not contain $L(w_0)$ as a 
subquotient either. Finally, we observe that $\mathrm{Rad}(N'(v_{i-1}))/\mathrm{Soc}(N'(v_{i-1}))$
contains $L(v_{i-1})$ as a subquotient while $\mathrm{Rad}(N'(w_{i+1}))/\mathrm{Soc}(N'(w_{i+1}))$
does not contain $L(v_{i-1})$ as a subquotient. This implies that $N'(v_{i-1})\not\cong N'(w_{i+1})$.
Hence, applying Lemma~\ref{hlem}, we obtain that the quiver of $Q(w)$ contains at 
least two loops at the point $w$.  This quiver also contains the elements $\alpha$
and $\beta$ described above. Factoring, if necessary, the extra arrows out, $Q(w)$ 
thus gives rise to the following configuration:
\begin{equation}\label{wconf1}
\xymatrix{
w_0\ar@{-}[rrd] && w\ar@{-}[lld]\ar@/^/@{-}[d]\ar@/_/@{-}[d]\\
w_0 && w
}.
\end{equation}
Since this is not an extended Dynkin quiver, this configuration is wild, see \cite{DF,DR0}.
Hence $\mathtt{D}(w)$, and thus $\HH$ is wild in this case.

{\bf Case~D.} Let $w=w_s$. In this case from \cite[\S~14]{ES} we have
$\mathrm{Ext}_{\mathcal{O}_{\lambda}}^1(L(w),L(v_{s-1}))\neq 0$. We also have
$\mathrm{Ext}_{\mathcal{O}_{\lambda}}^1(L(w),L(u_i))\neq 0$, $i=1,2$, since
$w_s$ and $u_i$ are neighbors. Hence the module $P(w)$ contains exactly $3$
copies of $L(w)$ in degree $2$: each lying in the top of the radical of some 
of the Verma modules $\Delta(x)$, $x=u_1,u_2,v_{s-1}$, occurring in degree $1$
in the Verma filtration of $P(w)$. Note that $L(w)$ does not occur in degree
$1$ (see Case~C). Further, $L(w_0)$ occurs at most one time in degree $1$ 
(this happens if $s=1$, in which case the occurrence in degree $1$ corresponds 
to the socle of $\Delta(w)$). In any case, since we have $3$ occurrences of $L(w)$ 
in degree $2$, at most one occurrence of $L(w_0)$ in degree $1$, and since 
$\mathrm{Ext}_{\mathcal{O}_{\lambda}}^1(L(w),L(w_0))\cong \mathbb{C}$ in the
case  $s=1$, mapping the degree $2$-occurrences to $I(w)$ we obtain at least
two non-isomorphic modules, $N_1$ and $N_2$, which have simple top and socle 
isomorphic to $L(w)$ and no other occurrences of $L(w)$ and $L(w_0)$. Taking
into account $\alpha$ and $\beta$, from Lemma~\ref{hlem} it now follows that 
some quotient of $Q(w)$ gives rise to the wild configuration \eqref{wconf1}. 
Hence $\mathtt{D}(w)$, and thus $\HH$ is wild in this case as well.

{\bf Case~E.} Finally, let $w=w_1$ and $s>1$. In this case both $\alpha$
and $\beta$ have degree $1$. From \cite[\S~14]{ES} we have
$\mathrm{Ext}_{\mathcal{O}_{\lambda}}^1(L(w),L(v_{0}))\neq 0$, which gives
us $2$ occurrences of $L(w)$ in degree $2$ of the module $P(w)$. One of them
comes from the subquotient $\Delta(v_0)$ in the Verma flag of $P(w)$. But
$v_0$ is dominant, and hence $\Delta(v_0)$ is in fact a submodule. Denote
by $\gamma$ the endomorphism of $P(w)$ of degree $2$, which corresponds to
this occurrence of $L(w)$ in $\Delta(v_0)$. Since
$(\beta\alpha)^2\neq 0$ by \cite[7.12-7.16]{GP}, it follows that 
the image of $\alpha\beta$ contains some $L(w_0)$ in degree $3$. However,
$\Delta(v_0)$ does not contain any $L(w_0)$ in degree $2$ (note that 
$\Delta(v_0)$ itself starts in degree $1$ in $P(w)$). Hence $\alpha\beta$ and $\gamma$
are linearly independent and thus $\gamma$ does not belong to the square of the radical. 
Now we claim that $\gamma^2=\gamma\alpha\beta=\alpha\beta\gamma=0$. The first and the
second equalities, that is $\gamma^2=\gamma\alpha\beta=0$, follow from the easy 
observation that $\Delta(v_0)$ does not have any $L(w)$ in degree $3=1+2$. 
The last one, that is $\alpha\beta\gamma=0$, follows from the fact that the
degree $1$-copy of $\Delta(v_0)$ belongs to the kernel of $\beta$ since
$P(w_0)$ does not have any $L(v_0)$ in degree $2$. Now, $P(w)$ has two copies of
$L(w)$ in the degree $2s$ which correspond to the subquotients $\Delta(u_1)$
and $\Delta(u_2)$ in the Verma flag of $P(w)$. Hence there should exist an
endomorphism of $P(w)$ of degree $2s$, which is linearly independent with $\alpha\beta$.
Since $\gamma^2=\gamma\alpha\beta=\alpha\beta\gamma=0$, it follows that 
this new endomorphism does not belong to the square of the radical of $Q(w)$.
Taking into account $\alpha$ and $\beta$, from Lemma~\ref{hlem} it now follows
that some quotient of $Q(w)$ gives rise to the wild configuration \eqref{wconf1}. 
Hence $\mathtt{D}(w)$, and thus $\HH$ is wild in this case as well.

This completes the cases \eqref{s3.n2}, \eqref{s3.n3}, \eqref{s3.n4}, \eqref{s3.n5}. 

Finally, let us consider the case \eqref{s3.n1}. Let $t_1$ and $t_2$ be the simple reflections
in $\mathbf{W}$, and let $\theta_{t_1}$, $\theta_{t_2}$ be translation functors through the
$t_1$ and $t_2$-wall respectively. If $\mathbf{H}\neq \mathbf{W}$,
then $\HH$ necessarily contains an indecomposable projective module,
which corresponds to some $w$ such that $\mathfrak{l}(w_0)-\mathfrak{l}(w)=2$.
The modules $\theta_{t_1}L(w_0)$ and $\theta_{t_2}L(w_0)$ are indecomposable and
have the following Loewy filtrations:
\begin{displaymath}
\begin{array}{lc}
 & L(w_0)\\
\theta_{t_1}L(w_0):&L({t_1}'w_0)\\
 & L(w_0)\\
\end{array},
\quad\quad
 \begin{array}{lc}
 & L(w_0)\\
\theta_{t_2}L(w_0):&L({t_2}'w_0)\\
 & L(w_0)\\
\end{array},
\end{displaymath}
for some ${t_1}',{t_2}'$ such that $\{t_1,t_2\}=\{t_1',t_2'\}$ (the exact values of 
$t_1'$ and $t_2'$ depend on the type of $\mathbf{W}$). 
In particular, $\theta_{t_1}L(w_0)\not\cong\theta_{t_2}L(w_0)$, both have simple top and simple 
socle isomorphic to $L(w_0)$, and both do not contain any subquotient isomorphic to 
$L(w)$ since $\mathfrak{l}(w_0)-\mathfrak{l}(w)=2$. 
Hence from Lemma~\ref{hlem} it  follows that the quotient of the corresponding $\mathtt{D}(w)$ 
modulo the square of the radical gives rise to the wild configuration \eqref{wconf2}. 
Hence $\mathtt{D}(w)$ is wild in this case. This proves Theorem~\ref{tmain}\eqref{tm.2}.

To complete the proof we just note that Theorem~\ref{tmain}\eqref{tm.3} follows from
Theorem~\ref{tmain}\eqref{tm.1} and Theorem~\ref{tmain}\eqref{tm.2} using the Tame
and Wild Theorem.

\section{The case of a semi-simple algebra $\mathfrak{g}$}\label{s4}

Theorem~\ref{tmain} is formulated for a simple algebra $\mathfrak{g}$. However, in the 
case of a semi-simple algebra the result is almost the same. In a standard way it
reduces to the description of the representation types of the tensor products of
algebras, described in Theorem~\ref{tmain}.

\begin{theorem}\label{tss}
Let $k>1$ be a positive integer, and $\mathtt{X}_i$, $i=1,\dots,k$, be basic algebras 
associated to non-semi-simple categories from the list of Theorem~\ref{tmain}. Then 
the algebra $\mathtt{X}_1\otimes\dots\otimes\mathtt{X}_k$ is never of finite representation 
type, and it is of tame representation type only in the following two cases:
\begin{enumerate}[(1)]
\item\label{tss.1} $k=2$ and both $\mathtt{X}_1$ and $\mathtt{X}_2$ have Coxeter type
$(A_1,e,A_1)$;
\item\label{tss.2} $k=2$, one of $\mathtt{X}_1$ and $\mathtt{X}_2$ has Coxeter type
$(A_1,e,A_1)$, and the other one has Coxeter type $(A_1,e,e)$.
\end{enumerate}
\end{theorem}

\begin{proof}
The algebra in \eqref{tss.1} is isomorphic to $\mathbb{C}[x,y]/(x^2,y^2)$ and hence
is tame with well-known representations. Let us thus consider the algebra $\mathtt{X}$
of the case \eqref{tss.2}. This algebra is given by the following quiver with relations
\begin{equation}\label{eqbre}
\xymatrix{  1 \ar@(ul,dl)[]_x  \ar@/^/[rr]^u   && 2 \ar@/^/[ll]^v  \ar@(ur,dr)[]^y }\qquad  x^2=y^2=uv=0,\ ux=yu,\ xv=vy.
\end{equation}

\begin{lemma}\label{lbre}
The algebra of \eqref{eqbre} is tame.
\end{lemma}

\begin{proof}
This algebra is tame by \cite{Be}, however, since the last paper is not easily available 
and does not contain a complete argument, we prove the tameness of $\mathtt{X}$. Consider 
the subalgebra $\mathtt{X}'\subset \mathtt{X}$ generated by $x,y,u$. Its indecomposable 
representations are 
\begin{align}\label{xyu}
\xymatrix@R=2ex{ {e_8} \ar[d] \\ {e_1} } \hskip1.5cm 
\xymatrix@R=2ex{ {e_9}\ar[d]\ar[dr] & {f_{10}}\ar[d] \\ {e_2} &{f_3} }\hskip1.5cm
\xymatrix@R=2ex{ {e_{10}}\ar[d]\ar[dr] & \\ {e_3}&{f_6} } \hskip1.5cm
\xymatrix@R=2ex{ {e_{11}}\ar[d]\ar[r] &{f_8}\ar[d] \\ {e_4}\ar[r] &{f_1} } \notag
\\
\\
\xymatrix@R=2ex{ {e_5} \\ {} } \hskip1.5cm 
\xymatrix@R=2ex{ {e_6} \ar[dr] & {f_9}\ar[d] \\ & {f_2} } \hskip1.5cm
\xymatrix@R=2ex{ {e_7}\ar[r] & {f_5} } \hskip1.5cm
\xymatrix@R=2ex{ {f_{11}}\ar[d] \\ {f_4} } \hskip1.5cm
\xymatrix@R=2ex{ {f_7} \\{} } \notag
\end{align} 
Here the elements $e_i$ form a basis of the space corresponding to the vertex $1$, the
elements $f_j$ form a basis of the space corresponding to the vertex $2$, the vertical
arrows show the action of $x$ and $y$, and the arrows going from left to right
show the action of $u$. Let $M$ be an $\mathtt{X}$-module. Decompose it as $\mathtt{X}'$-module.
Then the matrix $V$ describing the action of $v$ divides into the blocks $V_{ij},\
i,j=1,2,\dots,11$, corresponding to the basic elements $e_i$ and $f_j$ from above. 
Moreover, since $uv=0$, the blocks $V_{ij}$ can only be nonzero if $i\in\{1,2,3,5,8\}$;
since $xv=vy$, $V_{ij}=0$ if $i>4,j<5$ or $i>7,j<8$, and $V_{ij}=V_{i+7,j+7}$ for
$i,j\in\{1,2,3,4\}$. If $M'$ is another $\mathtt{X}$-module, $V'=(V'_{ij})$ is the corresponding
block matrix, a homomorphism $M\to M'$ is given by a pair of matrices $S,T$, where
$S:M(1)\to M(1),\,T:M(2)\to M(2)$. Divide them into blocks corresponding to the
division of $V$: $S=(S_{ij}),\,T=(T_{ij}),\ i,j=1,2,\dots,11$. One can easily check that
such block matrices define a homomorphism $M\to M'$ if and only if
the following conditions hold: 
\begin{itemize}
\item  $S$ and $T$ are block triangular, i.e. $S_{ij}=0$ and $T_{ij}=0$ if $i>j$.
\item  $S_{ij}=S_{i+7,j+7}$ and $T_{ij}=T_{i+7,j+7}$ for $i,j\in\{1,2,3,4\}$.
\item  $S_{ii}=T_{jj}$ if in the list \eqref{xyu} there is an arrow $e_i\to f_j$.
\item  $S_{ij}=T_{kl}$  if in the list \eqref{xyu} there are arrows $e_i\to f_k$ and $e_j\to f_l$. 
\item  $S_{ij}=0$ if $(i,j)\in\{(4,5),(4,6),(6,8),(7,8),(7,9)\}$.
\item   $T_{ij}=0$ if $(i,j)\in\{(3,5),(4,5), (4,6),(7,8)\}$.
\end{itemize}
Certainly, $S,T$ define an isomorphism if and only if all diagonal blocks are invertible.
In particular, we can replace the part $V_1=(V_{11}\ V_{12}\ V_{13}\ V_{14})$ by
$S_1^{-1}V_1T_1$, where $S_1$ is any invertible matrix and $T_1=(T_{ij}),\ i,j\in\{1,2,3,4\}$
is any invertible block triangular matrix. So we can suppose that $V_1$ is of the form
\begin{displaymath}
\left( \begin{array}{cc|cc|cc|cc}
0&I^{(1)}&0&0&0&0&0&0\\
0&0&0&I^{(2)}&0&0&0&0\\
0&0&0&0&0&I^{(3)}&0&0\\
0&0&0&0&0&0&0&I^{(4)}\\
0&0&0&0&0&0&0&0
\end{array}  \right), 
\end{displaymath}
where the vertical lines show the division of $V_1$ into blocks, $I^{(k)}$ denote identity
matrices (of arbitrary sizes).  Denote the parts of the blocks $V_{1j}$ to the right of
$I^{(k)}$ by $V_{1k,j}$ and those to the right of the zero part of $V_1$ by $V_{5j}$.
Using automorphisms, we can make zero all $V_{11,j}$ and $V_{12,j}$, as well as the
blocks $V_{13,j}$ and $V_{14,j}$ for $j>6$. Note that $V_{1j}=V_{8,j+7}$,
and we can also make zero all parts of the blocks $V_{1,j+7}$ over the parts $I^{(j)}$
of the blocks $V_{8,j+7}$. Subdivide the blocks of $S$ and
$T$ corresponding to this subdivision of $V_1$. Note that, since $S_{22}=T_{99}=T_{33}$,
we must also subdivide the blocks $S_{2j}$ into $S_{20,j}$ and $S_{21,j}$
respective to the zero and nonzero parts of $V_{13}$. Then the extra conditions for the
new blocks are: 
\begin{displaymath}
S_{21,20}=0 \quad\text{and}\quad S_{1k,1l}=0\quad \text{if}\quad k>l.
\end{displaymath}
Therefore, we get a matrix problem considered in \cite{Bo}. It is described by the
semichain
\begin{displaymath} 
\xymatrix@R=1ex{ f_5\ar[r] & f_6\ar[r]\ar[dr] & f_7\ar[dr] \\
&& f_8 \ar[r] & f_9\ar[r] & f_{10} \ar[r]& f_{11} }
\end{displaymath} 
for the columns, the chain
\begin{displaymath} 
e_5 \to  e_3 \to e_{21} \to e_{20} \to e_{15} \to e_{14} \to e_{13}
\end{displaymath} 
for the rows, and the unique equivalence $e_3\sim f_6$. This matrix problem is tame, hence,
the algebra $\mathtt{X}$ is tame as well.
\end{proof}

If $k>2$ then each of $\mathtt{X}_1$, $\mathtt{X}_2$, and $\mathtt{X}_3$ has 
at least one projective module with non-trivial endomorphism ring and thus 
$\mathtt{X}_1\otimes \mathtt{X}_2\otimes \mathtt{X}_3$ contains a centralizer subalgebra,
which surjects onto $\mathbb{C}[x,y,z]/(x,y,z)^2$. The later algebra is wild
by \cite{Dr2} and hence $\mathtt{X}$ is wild. 

If $k=2$ but none of the conditions \eqref{tss.1}, \eqref{tss.2} is satisfied, 
then one of the algebras $\mathtt{X}_1$ and $\mathtt{X}_2$ has a projective 
module, whose endomorphism algebra surjects onto $\mathbb{C}[x]/(x^3)$, and the other
one has a projective module, whose endomorphism algebra surjects onto $\mathbb{C}[y]/(y^3)$.
Hence there is a centralizer subalgebra in $\mathtt{X}$, which surjects onto
$\mathbb{C}[x,y]/(x^3,y^2)$, the later being wild by \cite{Dr2}. This shows that 
$\mathtt{X}_1\otimes \mathtt{X}_2$ is wild as well and completes the proof.
\end{proof}

\vspace{1cm}

\begin{center}
{\bf Acknowledgments}
\end{center}

The research was done during the visit of the first author to Uppsala
University, which was partially supported by the Faculty of Natural Science, 
Uppsala University, the Royal Swedish Academy of Sciences, and  The Swedish 
Foundation for International Cooperation in Research and Higher Education 
(STINT). These supports and the hospitality of Uppsala University are gratefully 
acknowledged. The second author was also partially supported by the Swedish 
Research Council. We would like to thank Catharina Stroppel for her comments
on the paper. We are especially in debt to the referee for a very careful
reading of the manuscript and for many useful comments, suggestions, and 
corrections.

\vspace{0.5cm}

\noindent
Yuriy Drozd, Department of Mechanics and Mathematics, Kyiv Taras
Shevchenko University, 64, Volodymyrska st., 01033, Kyiv, Ukraine,
e-mail: {\tt yuriy\symbol{64}drozd.org},\\
url: {\tt http://bearlair.drozd.org/$\sim$yuriy}.
\vspace{0.3cm}

\noindent
Volodymyr Mazorchuk, Department of Mathematics, Uppsala University,
Box 480, 751 06, Uppsala, SWEDEN, 
e-mail: {\tt mazor\symbol{64}math.uu.se},
url: {\tt http://www.math.uu.se/$\sim$mazor}.

\end{document}